\documentclass[12pt]{amsart}

\usepackage{amssymb}
\usepackage{amsmath}
\evensidemargin\oddsidemargin

\begin{document}
\setcounter{page}{1}

\newtheorem{PROP}{Proposition}
\newtheorem{REMS}{Remark}
\newtheorem{LEM}{Lemma}
\newtheorem{THE}{Theorem\!\!}
\newtheorem{COR}{Corollary}
\newtheorem{DEF}{Definition\!\!}
\newtheorem{CONJ}{Conjecture\!\!}

\renewcommand{\theTHE}{}
\renewcommand{\theCONJ}{}
\renewcommand{\theDEF}{}

\newcommand{\eqnsection}{
\renewcommand{\theequation}{\thesection.\arabic{equation}}
    \makeatletter
    \csname  @addtoreset\endcsname{equation}{section}
    \makeatother}
\eqnsection

\def\a{\alpha}
\def\b{\beta}
\def\B{{\bf B}} 
\def\BB{{\mathcal{B}}} 
\def\C{{\bf C}} 
\def\GG{{\mathcal{G}}} 
\def\II{{\mathcal{I}}}
\def\KK{{\mathcal{K}}} 
\def\LL{{\mathcal{L}}} 
\def\SS{{\mathcal{S}}}
\def\UU{{\mathcal{U}}}
\def\ca{c_{\a}}
\def\ka{\kappa_{\a}}
\def\coa{c_{\a, 0}}
\def\cua{c_{\a, u}}
\def\cL{{\mathcal{L}}} 
\def\Ea{E_\a}
\def\TT{{\bf T}}
\def\eps{{\varepsilon}} 
\def\esp{{\mathbb{E}}} 
\def\Ga{{\Gamma}} 
\def\G{{\bf \Gamma}} 
\def\e{{\rm e}}
\def\hr{\hat{\rho}} 
\def\ii{{\rm i}}
\def\L{{\bf L}}
\def\lbd{\lambda}
\def\lcr{\left[}
\def\lpa{\left(}
\def\lva{\left|}
\def\M{{\bf M}}
\def\MM{{\mathcal M}}
\def\NN{{\mathbb{N}}} 
\def\pb{{\mathbb{P}}}
\def\rl{{\mathbb{R}}}
\def\rpa{\right)}
\def\rcr{\right]}
\def\rva{\right|}
\def\W{{\bf W}}
\def\X{{\bf X}}
\def\XX{{\mathcal X}}
\def\YY{{\mathcal Y}}
\def\U{{\bf U}}
\def\V{{\bf V}_\a}
\def\Un{{\bf 1}}
\def\Z{{\bf Z}}
\def\A{{\bf A}}
\def\AA{{\mathcal A}}
\def\hAA{{\hat \AA}}
\def\hL{{\hat L}}
\def\hT{{\hat T}}

\def\claw{\stackrel{d}{\longrightarrow}}
\def\elaw{\stackrel{d}{=}}
\def\slaw{\stackrel{d}{\simeq}}
\def\qed{\hfill$\square$}

\title[On the law of homogeneous stable functionals]
      {On the law of homogeneous stable functionals}

\author[Julien Letemplier]{Julien Letemplier}

\address{Centre International de Valbonne, 190 Rue Fr\'ed\'eric Mistral, 06560 Valbonne, France. {\em Email} : {\tt ju.letemplier@gmail.com}}

\author[Thomas Simon]{Thomas Simon}

\address{Laboratoire Paul Painlev\'e, Universit\'e de Lille 1, Cit\'e Scientifique, 59655 Villeneuve d'Ascq, France. {\em Email} : {\tt simon@math.univ-lille1.fr}}



\begin{abstract} Let $\AA$ be the $\LL^q-$functional of a stable L\'evy process starting from one and killed when crossing zero. We observe that $\AA$ can be represented as the independent quotient of two infinite products of renormalized Beta random variables. The proof relies on Markovian time change, the Lamperti transform, and an explicit computation on perpetuities of hypergeometric L\'evy processes previously obtained by Kuznetsov and Pardo. This representation allows to retrieve several factorizations previously obtained by various authors, and also to derive new ones. We emphasize the connections between $\AA$ and more standard positive random variables. We also investigate the law of Riemannian integrals of stable subordinators. Finally, we derive several distributional properties of $\AA$ related to infinite divisibility, self-decomposability, and the generalized Gamma convolutions. 
\end{abstract}

\maketitle

\section{Introduction and statement of the results}

Let $L = \{L_t, \, t\ge 0\}$ be a real strictly $\a-$stable L\'evy process starting from one, and having characteristic exponent
\begin{equation}
\label{Norm}
\Psi(\lbd)\; =\;\log(\esp[e^{\ii \lbd L_1}])\; =\;\ii \lbd\; -\; (\ii \lbd)^\a e^{-\ii\pi\a\rho\, {\rm sgn}(\lbd)}, \qquad \lbd\in\rl,
\end{equation}
where $\a\in (0,2]$ is the self-similarity parameter and $\rho = \pb[L_1 -1\ge 0]$ is the positivity parameter. Recall that when $\a = 2,$ one has $\rho = 1/2$ and $\Psi(\lbd) = \ii \lbd -\lbd^2,$ so that $L = \sqrt{2} B$ is a rescaled Brownian motion. When $\a = 1,$ one has $\rho\in (0,1)$ and $L$ is a Cauchy process with a linear drift. When $\a\in (0,1)\cup(1,2)$ the characteristic exponent reads
$$\Psi(\lbd) \; =\; \ii \lbd\; -\; \kappa_{\a,\rho}\vert\lbd\vert^\a (1 - \ii\b\tan(\pi\a/2)\,{\rm sgn}(\lbd)),$$
where $\b\in[-1,1]$ is an asymmetry parameter, whose connection with the positivity parameter is given by Zolotarev's formula:
$$\rho \; =\; \frac{1}{2} \,+ \,\frac{1}{\pi\a} \arctan(\b\tan(\pi\a/2)),$$
and $\kappa_{\a,\rho} = \cos(\pi\a(\rho -1/2)) > 0$ is a scaling constant. We refer to Chapter 3 in \cite{S} for more details on this normalization, which could be modified without incidence on our purposes below. One has $\rho \in [0,1]$ if $\a < 1$ and $\rho\in[1-1/\a, 1/\a]$ if $\a > 1.$ When $\a > 1$ and $\rho = 1/\a$ the process $L$ has no positive jumps, whereas it has no negative jumps when $\a > 1$ and $\rho = 1-1/\a$. When $\a < 1,$ the process $L$ is a stable subordinator if $\rho =1$ and the opposite of a subordinator if $\rho =0.$ It is well-known that the first passage time below zero
$$T\; =\; T(\a,\rho)\; =\; \inf\{ t>0, \, L_t <0\}$$ 
is always a.s. finite, except in the subordinator case $\{\a <1, \rho = 1\},$ where it is a.s. infinite Consider now the following homogeneous functionals, possibly taking infinite values
$$\AA\;=\;\AA(\a,\rho,q)\; =\; \int_0^T \vert L_s\vert^q\, ds\; =\; \int_0^T (L_s)^q\, ds,\qquad q\in \rl.$$
In the case $q=0,$ one simply has $\AA =T$ and the law of this random variable has been the object of numerous works, old and recent - see \cite{Da, Bi, K1, K2, GJ, KP}, among others. It turns out that except in the simple case when $L$ is spectrally positive and $T$ is a positive stable random variable, the law of $T$ is quite complicated. In general, the random variables $\AA(\a,\rho,q)$ interpolate between the cases $q = -\infty$ and $q=+\infty,$ where a standard approximation shows that 
$$\AA(\a, \rho,q)^\frac{1}{q} \to \inf\{Z_t,\, t<T\}, \; q \to -\infty\quad \text{and}\quad \AA(\a, \rho,q)^\frac{1}{q} \to \sup\{Z_t,\, t<T\},\; q \to +\infty.$$
The law of these two stopped extrema, which we will respectively denote by $\II(\a,\rho)$ and $\SS (\a, \rho),$ is in contrast very simple. In the case $\a =2$, we can appeal to a well-known result on positive continuous martingales tending to zero to deduce that $\SS(2,1/2)^{-1}$ has a uniform distribution on $(0,1).$ In the case $\a<2$, the law of $\SS(\a, \rho)$ can be obtained from the results of Rogozin \cite{Ro} on the two-sided exit problem for stable L\'evy processes: one finds
\begin{equation}
\label{Rog}
\SS(\a,\rho)\; \elaw\; \B_{\a\hr, \a\rho}^{-1}
\end{equation}  
where, here and throughout, we set $\hr = 1-\rho$ and $\B_{a,b}$ denotes a Beta random variable with density
$$\frac{\Ga(a+b)}{\Ga(a)\Ga(b)} \, x^{a-1}(1-x)^{b-1} \, \Un_{(0,1)}(x).$$
In the remainder of this paper, we will use the conventions $\B_{0,b} \equiv 0, \B_{a,0} \equiv 1,$ and $\B_{1,1} =\U$ for the uniform random variable on $(0,1)$. To obtain the law of the stopped infimum in the non-trivial case with negative jumps, one may apply a computation performed by Port - see \cite{Po} or Exercise VIII.3 in \cite{Be1} - on the harmonic measure of the half-line, and deduce after some elementary manipulations the identity
\begin{equation}
\label{Port}
\II(\a,\rho)\;\elaw\;\B_{1-\a\hr, \a\hr}.
\end{equation}
Observe that our conventions allow to include the spectrally positive case $\{\a > 1, \hr = 1/\a\},$ where $\II(\a, 1-1/\a) \equiv 0,$ and the subordinator case $\{\a < 1, \rho =1\},$ where $\II(\a, 1) \equiv 1.$ In the case with negative jumps we also recover the well-known fact, following e.g. from the more general Lemma VIII.1 and Theorem VI.19 (ii) in \cite{Be1}, that $L$ must cross zero by a jump. This readily entails that $\AA$ is then a.s. finite for all $q\in\rl.$\\

In the spectrally positive case however, one has $L_T \equiv 0$ and the convergence of the integral defining $\AA$ is determined by the behaviour of the time-reversed process $\{L_{(T-t)-}, t < T\}$ at zero. By Lemma II.2 and Theorem VII.18 in \cite{Be1}, it is known that this process is distributed as the dual process ${\hat L} = -L$ conditioned to stay positive, starting from zero, until its last passage time above one. Moreover, the almost sure behaviour of this latter process at zero has been precisely described in \cite{Par} and it is easy to deduce from Theorem 1 and Theorem 3 i) therein (see also Theorem 2 and its proof in \cite{Be0}) that
$$\AA(\a, 1-1/\a,q)\; <\; +\infty \;\;\mbox{a.s.}\;\Longleftrightarrow\; \a + q > 0,$$
and that $\AA(\a, 1-1/\a,q) = +\infty$ a.s. if $\a + q \le 0.$ In the subordinator case, the laws of the iterated logarithm at infinity - see e.g. Theorems III.13-14 in \cite{Be1} - finally show that
$$\AA(\a, 1,q)\; <\; +\infty \;\;\mbox{a.s.}\;\Longleftrightarrow\; \a + q < 0,$$
and that $\AA(\a, 1,q) = +\infty$ a.s. if $\a + q \ge 0.$\\

Apart from the important case $q=0$ and the two limiting cases $q =\pm \infty,$ the random variable $\AA(\a, \rho,q)$ has also been investigated in the spectrally positive case with $q =-1$ in \cite{KPa}. The observations made in Section 3.1 of \cite{KPa} lead to the simple identity
\begin{equation}
\label{Fresh}
\AA(\a,1-1/\a, -1)\;\elaw\; \frac{1}{(\a-1) \L^{\a -1}}
\end{equation}
where, here and throughout, $\G_a$ denotes a Gamma random variable with density
$$\frac{x^{a-1} e^{-x}}{\Ga(a)} \, \Un_{(0,\infty)}(x),$$
and we have set $\L=\G_1$ for the standard exponential random variable. In our previous paper \cite{LS2}, using some connections with the stable Kolmogorov process, we have obtained a simple identity in law for $\AA(\a, 1-1/\a,1)$ in terms of a positive stable random variable and a Beta random variable. In this paper, we will give a general identity in law for all $\AA(\a,\rho,q)$'s in terms of Beta random variables only. We first  introduce the following definition. Throughout, it will be implicitly assumed that all products and quotients of given random variables are independent.

\begin{DEF} 
\label{defi}
For every $a,b,c > 0,$ the following a.s. convergent product
$$\TT(a,b,c)\;=\;\prod_{n=0}^\infty a_n\, \B_{a+nb,c}$$
with $a_n = (a+nb+c)/(a+nb),$ will be called an infinite Beta product with parameters $a,b,c.$ 
\end{DEF}

\bigskip

The almost sure convergence allowing $\TT(a,b,c)$ to be well-defined is a simple consequence of the martingale convergence theorem, and will be established in the next section. In the case $b > c,$ this is also the consequence of a more general resultat by Hackmann and Kuznetsov - see Proposition 1 in \cite{HK} - involving two increasing and unbounded sequences of parameters satisfying a certain interlacing property. As for Beta laws, we will adopt the conventions $\TT(0,b,c)\equiv 0$ and $\TT(a,b,0)\equiv 1.$  We can now state our main result.

\newpage

\begin{THE} 
\label{Main}
{\em (a)} If $L$ is spectrally positive, then for $q > -\a$ one has 
$$\AA(\a,1-1/\a,q)\; \elaw\; \frac{1}{(\a +q) \Ga(\a)\,\TT\lpa\frac{1}{\a +q}, \frac{1}{\a +q}, \frac{\a-1}{\a +q}\rpa}\cdot$$

\medskip

\noindent
{\em (b)} If $L$ is not spectrally positive, one has the following identities.

\bigskip

{\em (i)} For $q > -\a,$ 
$$\AA(\a,\rho,q)\; \elaw\; \lpa \frac{\Ga(1+q+\a\rho)\Ga(\a\hr)}{\Ga(1+\a+q)\Ga(\a)}\rpa\times\;\frac{\TT\lpa 1, \frac{1}{\a +q}, \frac{1-\a\hr}{\a +q}\rpa}{\TT\lpa\frac{\a\hr}{\a +q}, \frac{1}{\a +q}, \frac{\a\rho}{\a +q}\rpa}\cdot$$

\medskip

{\em (ii)} For $q =-\a,$ 
$$\AA(\a,\rho,-\a) \; \elaw\; \lpa \frac{\Ga(\a\hr) \Ga(1-\a\hr)}{\Ga(\a)}\rpa\times\; \L.$$

\medskip

{\em (iii)} For $q < -\a,$ 
$$\AA(\a,\rho,q)\; \elaw\; \lpa \frac{\Ga(1-q-\a\rho)\Ga(1-\a\hr)}{\vert \a+q\vert \Ga(\vert q\vert)}\rpa\times\;\frac{\B_{1, \frac{\a\hr}{\vert \a +q\vert}}\times\, \TT\lpa \frac{\vert q\vert}{\vert \a +q\vert}, \frac{1}{\vert \a +q\vert}, \frac{1-\a\rho}{\vert \a +q\vert}\rpa}{\TT\lpa\frac{1-\a\hr}{\vert \a +q\vert}, \frac{1}{\vert \a +q\vert}, \frac{\a\hr}{\vert \a +q\vert}\rpa}\cdot$$

\end{THE} 

\bigskip

Observe that except in the spectrally positive and the subordinator cases, our factorizations involve independent products of the type
$$\TT(a,b,c)\,\times\,\TT^{-1} (a',b',c'),$$
which will be called {\em infinite double Beta products} in the remainder of this paper. For certain values of the parameters, these double products take a simpler form and it is possible to recognize standard random variables. We will mention these simplifications all along the paper, especially in the last section. \\

Our main tools to obtain the above factorizations are Markovian time-change and the Lamperti transform introduced in \cite{La}, which imply that $\AA$ can always be viewed as the first-passage time at zero of some positive self-similar Markov process, and hence as the perpetuity of a certain L\'evy process. The latter has explicit L\'evy-Khintchine characteristics previously computed in \cite{CC}, and we can appeal to the formula for the Mellin transform of hypergeometric L\'evy perpetuities given by \cite{KP} in terms of the double Gamma function, which is that of an infinite double Beta product. \\

We now proceed to some remarkably simple factorizations which can be derived from the above main result. First, in the case $q=0,$ one can give a very quick proof of a observation made in \cite{Do} for stable L\'evy processes in duality. Let $\hL$ be the dual process of $L,$ that is the stable process with parameters $(\a, \hr)$ and the same normalization as above, starting from one, and let $\hT$ be the first time it crosses zero. 

\begin{COR}[Doney] 
\label{Doney}
Assume $\rho >0.$ Then, one has
$$\frac{T}{\hat{T}}\;\elaw\;\B_{\hr,\rho}^{-1} -\, 1.$$
\end{COR}

\medskip

The surprising point in the above identity is that the right-hand side does not depend on $\a.$ It is possible to obtain another factorization which does not depend on $\rho$ either. 

\begin{COR} 
\label{Dondon}
Assume $\a <1.$ Then, one has
$$\frac{\AA(\a,\rho,-1)}{\AA(\a,\hr,-1)}\;\elaw\;\U^{-1} -\, 1.$$
\end{COR}

\medskip

It is also interesting to replace the factorization obtained in part (b)-(iii) of the Theorem  in the context of stable subordinators starting from zero. Let $\sigma^{(\a)} = L -1$ with $\{\a <1, \rho =1\}$ be such a subordinator, and recall that with our normalization its Laplace transform reads 
$$\esp[e^{-\lambda \sigma^{(\a)}_t}]\; =\; e^{-t\lambda^\a},\qquad t,\lambda \ge 0.$$
As mentioned before, the laws of iterated logarithm at infinity imply
$$\int_0^\infty \frac{dt}{{(1+\sigma^{(\a)}_t)}_{}^q}\; =\; +\infty\;\;\mbox{a.s.}\quad\Longleftrightarrow\quad q\le\a.$$
The next corollary gives the law of the above Riemann random integral, whenever it is finite.
\begin{COR} 
\label{Subutex}
For every $q > \a$ one has
$$\int_0^\infty \frac{dt}{{(1+\sigma^{(\a)}_t)}_{}^q}\;\elaw\;\lpa \frac{\Ga(q-\a)}{\Ga(q)}\rpa\times\; \TT\lpa \frac{q}{q-\a}, \frac{1}{q-\a}, \frac{1-\a}{q-\a}\rpa.$$
\end{COR}

\medskip

Taking the limit $q\to +\infty$, it is then possible to obtain the law of the perpetuity
$$\int_0^\infty e^{-\sigma^{(\a)}_t}\, dt,$$
whose logarithm is distributed as the marginal at time $1-\a$ (with non-explicit density) of the real L\'evy process whose marginal at time 1 is the Gumbel distribution. This fact had already been implicitly observed in \cite{BY1} - see below (\ref{Gumgum}) for details. \\

Our main result also entails the two following multiplicative factorizations, which can be viewed as counterparts to Corollaries \ref{Doney} or \ref{Dondon}, in a more general setting. 

\begin{COR} 
\label{dual}
{\em (a)} Assume $q > -\a.$ For every $\rho>\rho',$ there exists an explicit infinite double Beta product $\XX$ such that
$$\AA(\a,\rho,q)\; \elaw\; \XX\,\times\, \AA(\a,\rho',q).$$

{\em (b)} Suppose that $Z$ is not spectrally positive and that $q < -\a.$ For every $\rho <\rho',$ there exists an explicit infinite double Beta product ${\hat \XX}$ such that
$$\AA(\a,\rho,q)\; \elaw\; \B_{\frac{\vert \a\rho' +q\vert}{\vert \a +q\vert},\frac{\a(\rho'-\rho)}{\vert \a +q\vert}}\,\times\,  {\hat \XX}\,\times\, \AA(\a,\rho',q).$$

\end{COR}

We can also establish a factorization of the Wiener-Hopf type for the functional $\AA(\a,\rho,q),$ in the spirit of \cite{PS, PS1}. 

\begin{COR} 
\label{WH}
Assume that $\AA(\a,\rho,q)$ is finite a.s. There exists a spectrally negative L\'evy process with positive mean $\{Z_t, \, t\ge 0\}$ and a subordinator $\{\sigma_t, \, t\ge 0\}$ such that
$$\AA(\a,\rho,q)\; \elaw\; \lpa \int_0^\infty e^{-Z_t}\, dt\rpa \times \lpa \int_0^\infty e^{-\sigma_t}\, dt\rpa.$$
The L\'evy-Khintchine characteristics of the two processes $\{L_t, \, t\ge 0\}$ and $\{\sigma_t, \, t\ge 0\}$ are explicitly given in terms of the parameters $(\a,\rho,q).$
\end{COR}

We next display several distributional properties related to the random variable $\AA(\a, \rho,q),$ which is always assumed to be finite a.s. in the remainder of this section. Our first result deals with the question whether its logarithm is infinitely divisible (ID) or self-decomposable (SD).

\begin{COR}
\label{AD} 
The random variable $\log \AA(\a,\rho,q)$ is always {\em ID}. Moreover, it is {\em SD} whenever $\a(1+\hr)\ge (\a +q)\wedge 1,$ and it is not always {\em SD.}  
\end{COR}

In particular, this entails that $\log\AA(\a,\rho,q)$ is SD for all $\a\geq 1$ or $q\le 0,$ and in particular that $\log T$ is always SD, a worthy property showing their unimodality by Theorem 53.1 in \cite{S}. The unimodality of $\AA(\a,\rho,q)$ itself is an interesting open problem in general, especially in the case $q=0.$ In our previous paper \cite{LS1}, we have shown the unimodality of hitting times for $L$, in the relevant case $\a >1.$ It is easy to see from the Theorem that $\AA(\a,\rho,q)$ has always a density on $\rl^+,$ whenever it is finite a.s. The unimodality amounts to the fact that this density has a unique local maximum. The following corollary characterizes the particular case when the density is non-increasing. In the case $q=0,$ the sufficiency of our criterion  had been observed in Corollary 1.5 of \cite{PS}. 

\begin{COR} 
\label{Chinchin}
The density function of $\AA(\a,\rho,q)$ is non-increasing on $\rl^+$ if and only if $1-\a\hr\ge\a +q\ge 0,$ or $0\ge \a +q\ge -\a\hr,$ or $\{\rho = 0, \a +q \le 0\}.$ Moreover, under these conditions, the density of $\AA(\a,\rho,q)$ is always bounded at zero.
\end{COR}

Our last corollary deals with the infinite divisibility of $\AA(\a,\rho,q)$ itself, in the spectrally positive case $\{\a >1, \rho = 1-1/\a\}.$ Recall that a given distribution on $\rl^+$ is called a generalized Gamma convolution (GGC) if it is the weak and non-degenerate limit of a weighted sum of independent Gamma distributions. We refer to Chapters 3-5 in \cite{Bd} for a classic account on GGC distributions. Notice that all GGC distributions are self-decomposable and hence infinitely divisible. It is well-known - see e.g. Example 3.2.1 in \cite{Bd} - that the GGC property is fulfilled by the random variable
$$\AA(\a, 1-1/\a,0)\; \elaw\; \Z_{\frac{1}{\a}}$$ 
where, here and throughout, we have set $\Z_\mu$ for the positive $\mu-$stable random normalized such that
$$\esp[e^{-\lbd \Z_\mu}]\; =\; e^{-\lbd^\mu}, \qquad \lbd \ge 0, \,\mu\in (0,1).$$
From (\ref{Fresh}) and known results on inverse Gamma distributions - see Chapter 5 in \cite{Bd} and also Section 3.1 for a specific account, one can show that the law of $\AA(\a, 1-1/\a, -1)$ is a GGC. In our previous paper \cite{LS2} - see Corollary 2 therein, we had observed that the law of $\AA(\a, 1-1/\a, 1)$ is also a GGC. The following result generalizes these facts.

\begin{COR} 
\label{AGGC}
For every $\a > 1, q > -\a,$ the law of $\AA(\a,1-1/\a,q)$ is a {\em GGC}.
\end{COR}

It is interesting to notice that $\AA(\a,\rho,q)$ is not infinitely divisible in general. For every $\a< 1,$ the easily established identity in law
\begin{equation}
\label{MLL}
\AA(\a,0,0)\; \elaw\; \Z_\a^{-\a}
\end{equation}
shows indeed that $\A(\a,0,0)$ is not ID because it has superexponential tail distributions at infinity - see Exercices 29.18 and 29.19 in \cite{S}. Since infinite divisibility is preserved under weak convergence, for every $\a < 1$ we see that there exists $\rho > 0$ and $q\neq 0$ such that $\AA(\a,\rho,q)$ is not ID either. It seems difficult to characterize the infinite divisibility of $\AA(\a,\rho,q)$ in terms of the parameters $(\a,\rho,q).$ In the important case $q=0$, it is natural to raise the following conjecture.

\begin{CONJ} With the above notations, one has
$$T\;\mbox{{\em is ID}}\;\Longleftrightarrow\; \a \ge 1.$$
\end{CONJ}
The answer to this question does not seem immediate and we believe that the factorization obtained in Part (b)-(i) of the Theorem might help to solve the if part of it. This will be the matter of future research. \\

\section{Proof of the Theorem}

\subsection{Preliminary facts on infinite Beta products} 

Let us first check the a.s. convergence of the infinite product defining $\TT(a,b,c).$ To do so, we first observe that 
$${\tilde \TT} (a,b,c) \; =\; \prod_{n=0}^\infty {\tilde a}_n\, \B_{a+nb,c}$$
with ${\tilde a}_n = e^{\psi(a +nb+c) - \psi(a+nb)},$ where $\psi$ is the standard digamma function, is a.s. convergent. Indeed, a well-known consequence of Malmsten's formula for the Gamma function - see e.g. 1.9(1) p.21 in \cite{EMOT} for the latter formula - shows that
\begin{equation}
\label{IniBB}
\esp[\B_{a+nb,c}^s]\;=\;\exp[\int_{-\infty}^0 (e^{sx}-1)\frac{e^{-(a+nb)|x|}(1-e^{-c|x|})}{|x|(1-e^{-|x|})}\, dx]
\end{equation}
for every $a,b,c,s > 0$ and $n\in\NN.$ Differentiating this formula at $s =0$ and applying Gauss' integral formula 1.7.2(17) in \cite{EMOT} yields
$$\esp[\log \B_{a+nb,c}]\; =\; \psi(a+nb)\; -\; \psi(a+nb +c)\; =\; -\log a_n.$$
Differentiating a second time, we obtain
$${\rm Var}[\log \B_{a+nb,c}]\; =\; \int_0^{\infty} \frac{x(1-e^{-cx})}{(1-e^{-x})}\,e^{-(a+nb)x} dx.$$
All of this entails that the discrete martingale 
$$X_n\; =\; \sum_{k=1}^n \log (a_k \B_{a+kb,c}),\;\; n\ge 1,$$
converges a.s. because its variances are uniformly bounded by 
$$\int_0^{\infty} \frac{x(1-e^{-cx})}{(1-e^{-x})(1-e^{-bx})}\,e^{-ax} \, dx\; <\; +\infty.$$
Hence, the product ${\tilde \TT}(a,b,c) = \lim_{n\to+\infty} e^{X_n}$
converges a.s. Moreover, it is clear from the above considerations that its Mellin transform is given by
\begin{equation}
\label{MelTT}
\esp[{\tilde \TT}(a,b,c)^s]\;=\;\exp \int_{-\infty}^0 (e^{sx}-1 -sx)\frac{e^{-a|x|}(1-e^{-c|x|})}{|x|(1-e^{-|x|})(1-e^{-b\vert x\vert})}\, dx
\end{equation}
for every $s > -a.$ In particular, ${\tilde \TT}(a,b,c)$ is integrable and Fubini's theorem implies
$$\esp[{\tilde \TT}(a,b,c)]\; =\; \prod_{n=0}^\infty \lpa \frac{a +nb}{a+nb +c} \rpa e^{\psi(a +nb+c) - \psi(a+nb)}\; <\; \infty.$$
All of this shows that the infinite product
\begin{equation}
\label{Quo}
\TT(a,b,c)\; =\; \frac{{\tilde \TT}(a,b,c)}{\esp[{\tilde \TT}(a,b,c)]}
\end{equation}
is a.s. convergent.

\qed

\begin{REMS} {\em (a) Among all possible normalizing deterministic sequences making the infinite Beta product a.s. convergent, the above $\{ a_n\}$ and $\{{\tilde a}_n\}$ are such that 
$$\esp[\TT(a,b,c)]\; =\; 1\qquad\mbox{and}\qquad \esp[\log {\tilde \TT}(a,b,c)]\; =\; 0.$$
The random variable ${\tilde \TT}(a,b,c)$ is less simple than $\TT(a,b,c),$ but easier to handle for computations. From time to time, we will first prove an identity for ${\tilde \TT}(a,b,c)$ and then deduce the corresponding one for $\TT(a,b,c)$ by normalization. \\

(b) It is interesting to compare (\ref{MelTT}) with the Mellin transform of the Barnes Beta distributions recently introduced in \cite{O}. Notice that the latter have their support in $[0,1]$ whereas in (\ref{MelTT}), the integral of the L\'evy measure along $\vert x\vert \wedge 1$ is infinite, so that $\log \TT(a,b,c)$ has full support - see e.g.  Theorem 24.10 in \cite{S} -  and Supp $\TT(a,b,c)\, =\, \rl^+.$ }  
\end{REMS}

We next state some simple multiplicative identities in law for $\TT(a,b,c)$ which will be important in the sequel. We use standard notation for size-bias and scale mixtures. If $X$ is a positive random variable and $\nu\in\rl$ is such that $\esp[X^\nu] < \infty,$ the positive random variable $X^{(\nu)}$ is defined by 
$$\esp[f(X^{(\nu)})]\; =\; \frac{\esp[X^\nu f(X)]}{\esp[X^\nu]}$$
for all $f$ bounded continuous. If $X, Y$ are positive random variables such that $X\,\elaw\, c\, Y$ for some $c >0,$ we write $X\,\slaw\, Y.$ 

\begin{PROP}  One has
\label{GZT}
\begin{equation}
\label{MulTT}
\TT(a,b,c)\;\times\;\TT(a+c,b,d)\;\elaw\;\TT(a,b,c+d),
\end{equation}
\begin{equation}
\label{MulTTa}
\TT(a,b,c)\;\slaw\;\TT(ab^{-1},b^{-1},cb^{-1})^{\frac{1}{b}},
\end{equation}
\begin{equation}
\label{MulTTb}
\TT(a,b,c)\;\elaw\;\B_{a,c}\,\times\,\TT(a,b,c)^{(b)}\;\elaw\;\B_{\frac{a}{b},\frac{c}{b}}^{\frac{1}{b}}\,\times\,\TT(a,b,c)^{(1)},
\end{equation}
\begin{equation}
\label{MulTTc}
\G_a\;\elaw\; a\,\TT(a,b,b)\quad{\rm and}\quad \G_a^b\;\elaw\; \frac{\Ga (a+b)}{\Ga(a)}\,\TT(ab^{-1},b^{-1},1),
\end{equation}
\begin{equation}
\label{MulTTd}
\Z_\a^{-1}\;\elaw\;\Ga(1+\a^{-1})\,\TT(\a,\a,1-\a)\quad{\rm and}\quad \Z_\a^{-\a}\;\elaw\; \frac{1}{\Ga(1+\a)}\,\TT(1,\a^{-1},\a^{-1}-1).
\end{equation}

\end{PROP}

\proof The identity ${\tilde \TT}(a,b,c)\times{\tilde \TT}(a+c,b,d)\elaw{\tilde\TT}(a,b,c+d),$ which implies (\ref{MulTT}) by normalization, is a straightforward consequence of (\ref{MelTT}). Similarly, a change of variable inside (\ref{MelTT}) yields easily
$${\tilde \TT}(a,b,c)\;\elaw\;{\tilde \TT}(ab^{-1},b^{-1},cb^{-1})^{\frac{1}{b}},$$
which again implies (\ref{MulTTa}) by normalization. Since $\B_{a+ (n+1)b,c}\elaw\B_{a+ nb,c}^{(b)},$ the first identity in (\ref{MulTTb}) is plain, whereas the second one follows from (\ref{MulTTa}) and normalization. Finally, the two identities in (\ref{MulTTc}) and (\ref{MulTTd}) are consequences of Lemmas 2 and 3 in \cite{BS3},  (\ref{MulTTa}), and normalization.

\qed

\begin{REMS}{\em The identity (\ref{MulTT}) mimics the well-known and easily established identity
\begin{equation}
\label{MulBB}
\B_{a,b}\;\times\;\B_{a+b,c}\;\elaw\;\B_{a,b+c},
\end{equation}
which also be useful in the sequel.}
\end{REMS}

We finally compute the Mellin transform of $\TT(a,b,c)$ in a different fashion than putting (\ref{MelTT}) and (\ref{Quo}) together would do. The formula, which involves the double Gamma function $G(z, \tau)$ introduced in \cite{Ba0, Ba}, will make it possible to use the results of \cite{KP} in the non-spectrally positive case.

\begin{PROP}  
\label{Polvo}
For every $s > -a,$ one has
$$\esp[\TT(a,b,c)^s]\;=\; \lpa \frac{\Ga(ab^{-1})}{\Ga ((a+c)b^{-1})}\rpa^s\times\; \frac{G(a+c+s, b) G(a,b)}{G(a+s, b)G(a+c, b)}\cdot$$
\end{PROP}

\proof Setting 
$$\MM(s)\; = \; \lpa \frac{\Ga ((a+c)b^{-1})}{\Ga(ab^{-1})}\rpa^s\times\;\esp[\TT(a,b,c)^s],$$ 
we first observe that the second identity in (\ref{MulTTb}) implies the functional equation
\begin{equation}
\label{FE}
\MM(s+1)\; = \; \lpa \frac{\Ga ((a+c+s)b^{-1})}{\Ga((a+s)b^{-1})}\rpa\times\;\MM(s).
\end{equation}
On the one hand, the formula
\begin{equation}
\label{Gumbel}
\Gamma(\delta+s)\;=\; \exp[-\psi(\delta) s+ \int_{-\infty}^0 (e^{sx}-1-sx)\frac{e^{-\delta\vert x\vert}}{\vert x\vert (1-e^{\vert x\vert})}\, dx], \quad s > -\delta,
\end{equation}
shows that the function
$$s\;\mapsto\; \frac{\Ga ((a+c+s)b^{-1})}{\Ga((a+s)b^{-1})}$$
is log-concave on $(-a, \infty).$ Hence, by the main result of \cite{Wb}, we know that there is a unique log-convex solution to (\ref{FE}) satisfying $\MM(0) = 1.$ On the other hand, the concatenation formula $G(z+1, \tau) = \Ga(z\tau^{-1}) G(z, \tau)$ implies that 
$$s\;\mapsto\;\frac{G(a+c+s, b) G(a,b)}{G(a+s, b)G(a+c, b)}$$
is such a solution. This completes the proof, by the log-convexity of $\MM(s).$

\qed

\subsection{The Brownian case} In this classical and particular case, it is  possible to show the Theorem in two ways, an additive one and a multiplicative one. Recall that one has to  prove Part (a) only, with $\{\a = 2, \rho =1/2\}.$

\subsubsection{Additive proof} We will use the standard Feynman-Kac formula. Fix  $q >-2.$ From e.g. Formula (3.2) in \cite{JPY}, we know that for every $\lambda \ge 0,$ the quantity $\esp[e^{-\lambda \AA(2,1/2,q)}]$ is given by the value at $t=1$ of the unique positive non-increasing solution $u(t)$ to the following Sturm-Liouville problem on $\rl^+$:
\begin{equation}
\label{bracebrace}
\left\lbrace \begin{split} 
&u''(t)-\lambda \, t^q\,  u(t)\;=\;0,\\
&u(0)\;=\;1.
\end{split}
\right.
\end{equation}
The solution to this problem follows from known results on Bessel functions - see \cite{K}. Let us give the details for completeness. Setting $u(t)=v(\lambda^{\frac{1}{q+2}}t)$, we are reduced to
\[v''(t)\,-\, t^q\, v(t)\;=\;0.\]
Setting now $a=\frac{q+2}{2}$ and $v(t)=\sqrt{t}\,w(t^a)$, the equation becomes
$$a^2 \, t^2\, w''(t)\;+\;a^2\, t\, w'(t)\;-\;\lpa \frac{1}{4}+t^2\rpa w(t)\;=\;0,$$
whose positive non-increasing solutions on $\rl^+$ are well-known to be of the type $c\, K_{\frac{1}{q+2}}(a^{-1} t)$, where $K_\nu$ is the Macdonald function - see sections 7.2.2. and 7.12 in \cite{EMOT}. We deduce 
$$\esp\lpa e^{-\lambda \AA(2,1/2,q)}\rpa\;=\; c\, \lambda^{\frac{1}{2(q+2)}}K_{\frac{1}{q+2}}\lpa\frac{2 \sqrt{\lambda}}{q+2}\rpa,$$
where $c$ is the normalizing constant. Applying e.g. Exercice 34.13 in \cite{S}, we can invert the Laplace transform and find
\begin{equation}
\label{MainB}
\AA(2,1/2,q)\; \elaw \;\frac{1}{(q+2)^2\G_{\frac{1}{q+2}}}\cdot
\end{equation}
We can finally conclude thanks to the first identity in (\ref{MulTTc}). 

\qed

\begin{REMS} {\em (a) In the cases $q=1,2$ and in theses cases only, the Feynman-Kac method also allows to evaluate the bivariate Laplace transform of $(T, \AA(2, 1/2,q))$ - see Example 5 in the fourth section of \cite{JPY}. \\

(b) The Feynman-Kac method extends formally to stopped functionals of general L\'evy processes - see e.g. \cite{Be1} p. 140. However, it does not seem possible to solve explicitly the integro-differential equations appearing in the framework of our homogeneous stable functionals, when $\a < 2.$ }
\end{REMS}

\medskip

\subsubsection{Multiplicative proof} This proof is more elaborate than the additive one. However, it will be possible to extend it to the general stable case. Consider the killed process
$$X_t\; =\; L_t\, \Un_{\{T>t\}},$$
which is a continuous strong Markov process, non-negative, and self-similar with index $1/2$ in the sense of \cite{La}. For every $q > -2,$ the additive functional 
$$t\;\mapsto\;A_t\; =\; \int_0^t \vert X_s\vert^q\, ds$$
is continuous and increasing on $[0,T]$. Set $\vartheta_t = \inf\{ u > 0, \, A_u > t\}$ for its inverse, which is also continuous and increasing on $[0,\AA(2,1/2,q)].$ The time-changed process
$$Y_t \; = \; X_{\vartheta_t}$$
is continuous, non-negative, and also strong Markov - see \cite{Dy} p. 10 for the latter property. It is easy to see that $Y$ is also self-similar with index $1/(q+2).$ Since by construction one has $\AA (2,1/2,q) =\inf\{t>0,\, Y_t =0\},$ Theorem 4.1 in \cite{La} entails that
$$\AA (2,1/2,q) \; \elaw\; \int_0^\zeta e^{-(q+2)\,\xi_t}\, dt,$$
where $(\xi, \zeta)$ is a certain L\'evy process with lifetime $\zeta.$ To identify the latter, we observe from \cite{Dy} p. 13 that the action of the infinitesimal generator of $Y$ is deduced from that of $X$ by the formula
\begin{equation}
\label{ChangeB}
\LL^{}_Y f(x)\; =\; x^{-q}\,\LL^{}_X f(x) \; =\; x^{-q} f''(x)
\end{equation}
for each twice differentiable fonction $f : (0, +\infty)\to\rl.$ From Theorem 6.1 in \cite{La}, this entails that $(\xi, \zeta)$ is the Lamperti L\'evy process associated to $X$, which is known - see e.g. \cite{Yo} - to be the drifted Brownian motion $\{L_t + t, \, t \ge 0\}.$ We deduce that 
$$\AA (2,1/2,q) \; \elaw\; \int_0^\infty e^{-(q+2) (L_t + t)}\, dt\; \elaw\; \frac{1}{2(q+2)^2}\int_0^\infty e^{-(B_t + at/2)}\, dt$$
where $\{B_t, \, t \ge 0\}$ is a standard Brownian motion and $a = 1/(q+2)^2.$ Dufresne's identity on Brownian perpetuities (see \cite{Du}, and also \cite{BCGL} for a concomitant proof using arguments from theoretical physics) yields finally
$$\AA(2,1/2,q)\; \elaw \;\frac{1}{(q+2)^2\G_{\frac{1}{q+2}}}\cdot$$

\qed

\subsection{The spectrally positive case} Repeating the first part of the multiplicative Brownian proof, for every $q > -\a$ we show without difficulty that 
\begin{equation}
\label{ExpSP}
\AA (\a,1-1/\a,q) \; \elaw\; \int_0^\zeta e^{-(q+\a)\,\xi_t}\, dt,
\end{equation}
where $(\xi, \zeta)$ is the Lamperti L\'evy process associated with the positive self-similar Markov process 
$$X_t\; =\; L_t\, \Un_{\{T>t\}}.$$
Using Corollary 1 in \cite{CC} and the Lemma in \cite{BS1} - see also Lemma 1 in \cite{KPa}, we see that $\zeta =+\infty$ a.s. and that the Laplace exponent of $\xi$ is given by
$$\log \esp[e^{z \xi_1}]\; =\;  \frac{z\,\Ga(\a+z)}{\Ga(1+z)}$$
for every $z > -\a.$ Using the proof of Proposition 2 in \cite{BY2}, we deduce that the Mellin transform of $\AA (\a,1-1/\a,q)$ satisfies the functional equation
$$\esp[\AA(\a, 1-1/\a, q)^{-s-1}]\; =\;\frac{(\a+q)\Ga(\a+(\a+q)s)}{\Ga(1+(\a+q)s)}\,\times\,\esp[\AA(\a, 1-1/\a, q)^{-s}]$$ 
for every $s > -1/(\a+q).$ Hence, setting $\BB = ((\a+q)\Ga(\a) \AA(\a, 1-1/\a, q))^{-1},$ we obtain the multiplicative relation
$$\BB \; \elaw\; \B_{1, \a -1}^{\a+q}\; \times\; \BB^{(1)}.$$
By Theorem 3.5 in \cite{Pk}, this equation has a unique solution having unit expectation, and the second identity in (\ref{MulTTb}) shows that this solution must be
$$\TT\lpa\frac{1}{\a +q}, \frac{1}{\a +q}, \frac{\a-1}{\a +q}\rpa.$$
This completes the proof.

\qed

\begin{REMS} {\em (a) In the limiting case $\{\a =1, \rho =0\}$ where $Z_t = 1-t$ and $T\equiv 1,$ we find 
$$\AA(1, 0, q)\; \equiv\;\frac{1}{q+1}\; =\; \int_0^1 (1-t)^q\, dt$$ 
as expected. \\

(b) In the case $q =-1,$ the above result and the second identity in (\ref{MulTTc}) allow us to retrieve (\ref{Fresh}) readily. On the other hand, it is not straightforward to retrieve the main result of \cite{LS2} from the above factorization.}
\end{REMS}

We now proceed to Part (b), noticing first that we can exclude the case where $L$ is a subordinator. Indeed, one has $\AA(\a, 1,q) = +\infty$ if $q \ge -\a$ whereas in the case $q < -\a,$ the law of $\AA(\a, 1,q)$ is obtained from that of $\AA(\a, \rho,q)$ at the limit $\rho\uparrow 1.$ This is given as Corollary \ref{Subutex} below. 

\subsection{The case with negative jumps}

\subsubsection{The case $q>-\a$} Using the same argument as in the spectrally positive case combined with Corollary 1 in \cite{CC} and the full strength of Theorem 1 in \cite{KP} - which is stated in the spectrally two-sided case only, but readily applies to the two limiting cases $\{\a <1, \rho =0\}$ and $\{\a > 1, \rho =1/\a\}$, we first obtain the identity 
\begin{equation}
\label{Exp}
\AA(\a,\rho,q) \; \elaw\; \int_0^\zeta e^{-(q+\a)\,\xi_t}\, dt,
\end{equation}
where $\xi$ is a L\'evy process having Laplace exponent
$$-\log[\esp[e^{z \xi_1}]]\; =\;\frac{\Ga(1-z)\Ga(\a+z)}{\Ga(1-\a\hr-z)\Ga(\a\hr+z)}$$
for every $z\in (-\a,1),$ and whose lifetime $\zeta$ is here a.s. finite. More precisely, one has $\zeta \elaw r \L$ with
$$r \; =\; \frac{\Ga(1-\a\hr)\Ga(\a\hr)}{\Ga(\a)}\cdot$$
Theorem 2 in \cite{KP} entails then
\begin{equation}
\label{Frac}
\esp[\AA(\a,\rho,q)^s]\; =\; \Psi_-(s)\,\times\, \Psi_+(s)
\end{equation}
for every $s \in (-1, \frac{\a\hr}{\a+q})$, with 
$$\Psi_-(s)\; = \;\frac{G(\frac{\a}{\a +q}-s, \frac{1}{\a+q})G(\frac{\a\hr}{\a +q}, \frac{1}{\a+q})}{G(\frac{\a}{\a +q}, \frac{1}{\a+q})G(\frac{\a\hr}{\a +q}-s, \frac{1}{\a+q})}$$
and
$$ \Psi_+(s)\; =\; \Ga(1+s)\,\frac{G(\frac{1+\a\rho +q}{\a +q} +s, \frac{1}{\a+q})G(\frac{1+\a +q}{\a +q}, \frac{1}{\a+q})}{G(\frac{1+\a\rho +q}{\a +q}, \frac{1}{\a+q})G(\frac{1+\a +q}{\a +q} +s, \frac{1}{\a+q})}\cdot$$
On the one hand, Proposition \ref{Polvo} entails
$$\Psi_-(s)\; = \; \lpa \frac{\Ga(\a\hr)}{\Ga(\a)}\rpa^s\times\;\esp\lcr \TT\lpa\frac{\a\hr}{\a +q}, \frac{1}{\a +q}, \frac{\a\rho}{\a +q}\rpa^{-s}\rcr.$$ 
To handle the factor $\Psi_+(s),$ we appeal to a formula of the Malmsten type for Barnes' double Gamma function which was obtained in \cite{LK}, Formula (A.15), and which reads
\begin{equation}
\label{Malm}
G(s,\delta)\; =\; \exp\int_0^\infty \lpa \frac{1-e^{-sx}}{(1-e^{-x})(1-e^{-\delta x})} \, -\, \frac{se^{-\delta x}}{1-e^{-\delta x}}\, +\, (s-1)(\frac{s}{2\delta} -1)e^{-\delta x}\rpa\frac{dx}{x}
\end{equation}
for all $s,\delta > 0.$ From this formula, it is easy to see that for all $a,b,\delta > 0$ and $s < \inf(a,b),$ one has
\begin{equation}
\label{DoubleG}
\log\lpa\frac{G(b-s, \delta)G(a, \delta)}{G(b, \delta)G(a-s, \delta)}\rpa\; =\; \kappa\, s\; +\; \int_0^{\infty} (e^{sx} - 1 - sx) f_{a,b,\delta}(x) dx,
\end{equation}
for some normalizing constant $\kappa,$ and with the notation
$$f_{a,b,\delta}(x) \; =\; \frac{(e^{-ax}-e^{-bx})}{x(1-e^{-x})(1-e^{-\delta x})}\cdot$$
Combining (\ref{DoubleG}) and (\ref{Gumbel}), we deduce after some simplifications
$$\log \Psi_+(s)\; =\; (\kappa -\gamma) s \; +\; \int_{-\infty}^0 (e^{sx} - 1 -sx)\, f_{1,\frac{1+\a\rho+q}{\a+q},\frac{1}{\a+q}}(\vert x\vert)\, dx$$
for every $s>-1,$ with 
$$ a\; =\; \frac{1+\a\rho+q}{\a+q},\quad b\; =\; \frac{1+\a+q}{\a+q}\quad \mbox{et}\quad \delta\; =\; \frac{1}{\a+q}\cdot$$
Using again (\ref{MelTT}) and putting everything together lead finally to
$$\AA(\a,\rho,q)\; \elaw\; \KK\;\TT\lpa 1, \frac{1}{\a +q}, \frac{1-\a\hr}{\a +q}\rpa\;\times\;\TT^{-1}\lpa\frac{\a\hr}{\a +q}, \frac{1}{\a +q}, \frac{\a\rho}{\a +q}\rpa$$
with
$$\KK\; =\; \frac{\Ga(\a\hr)}{\Ga(\a)}\,\times\,\Psi_+(1)\; =\; \frac{\Ga(1+q+\a\rho)\Ga(\a\hr)}{\Ga(1+\a+q)\Ga(\a)}\cdot
$$

\qed

\subsubsection{The case $q=-\a$} This is a consequence of the case $\{q<-\a\}$ and continuity in law. Indeed, since $\zeta$ is a.s. finite, it is clear that the application
$$q\;\mapsto\; \int_0^\zeta e^{-(q+\a)\,\xi_t}\, dt$$
is a.s. continuous and hence continuous in law on $\rl$, with the above notation for $\xi.$ Since $q\mapsto \AA(\a,\rho, q)$ is also plainly continuous in law we obtain, letting $q\to -\a,$
$$\AA(\a,\rho, -\a)\;\elaw\; \zeta\; \elaw\; \lpa\frac{\Ga(1-\a\hr)\Ga(\a\hr)}{\Ga(\a)}\rpa\L.$$

\qed

\subsubsection{The case $q<-\a$} Consider the process
$${\hat X}_t \; =\; L^{-1}_t\, \Un_{\{T>t\}},$$
which is a positive self-similar Markov process having negative self-similarity parameter ${\hat \a} = -\a,$ and whose lifetime is
$$T\; =\; \int_0^\zeta e^{-\a\,\xi_t}\, dt\; =\; \int_0^{\hat \zeta} e^{-{\hat \a}\,{\hat \xi}_t}\, dt,$$
with the above notation for the L\'evy Lamperti process $(\xi,\zeta)$ and where $({\hat \xi},{\hat \zeta}) = (-\xi,\zeta)$ stands for its dual process. Introducing the homogeneous additive functional 
$$t\;\mapsto\;{\hat A}_t\; =\; \int_0^t \vert {\hat X}_s\vert^{\vert q\vert}\, ds$$
and using exactly the same time-change argument as above shows the identity
\begin{equation}
\label{Exp-}
\AA(\a,\rho,q) \; \elaw\; \int_0^{\hat \zeta} e^{-(\vert q\vert +{\hat \a})\,{\hat \xi}_t}\, dt\; =\; \int_0^\zeta e^{-\vert q+\a\vert \,{\hat \xi}_t}\, dt.
\end{equation}
On the other hand, the Laplace exponent of ${\hat \xi}$ is given by
$$-\log[\esp[e^{z {\hat \xi}_1}]]\; =\;\frac{\Ga(1+z)\Ga(\a-z)}{\Ga(1-\a\hr+z)\Ga(\a\hr-z)}$$
for every $z\in (-1, \a).$ The remainder of the proof is now analogous to the above, but the details are a bit different. Using again Theorem 2 in \cite{KP}, we obtain
\begin{equation}
\label{Frac-}
\esp[\AA(\a,\rho,q)^s]\; =\; \Psi_-(s)\,\times\, \Psi_+(s)
\end{equation}
for every $s \in (-1, \frac{1-\a\hr}{\vert \a+q\vert})$, where
\begin{eqnarray*}
\Psi_-(s) & = &\frac{G(\frac{1}{\vert \a +q\vert}-s, \frac{1}{\vert \a+q\vert})G(\frac{1-\a\hr}{\vert \a +q\vert}, \frac{1}{\vert \a+q\vert})}{G(\frac{1}{\vert \a +q\vert}, \frac{1}{\vert \a+q\vert })G(\frac{1-\a\hr}{\vert \a +q\vert}-s, \frac{1}{\vert \a+q\vert})}\\
& = &  \lpa \Ga(1-\a\hr)\rpa^s \esp\lcr \TT\lpa\frac{1-\a\hr}{\vert \a +q\vert}, \frac{1}{\vert \a +q\vert}, \frac{\a\hr}{\vert \a +q\vert}\rpa^{-s}\rcr
\end{eqnarray*}
and
$$ \Psi_+(s)\; =\; \Ga(1+s)\,\frac{G(\frac{\vert \a\rho +q\vert}{\vert\a +q\vert} +s, \frac{1}{\vert\a+q\vert})G(\frac{\vert q\vert}{\vert\a +q\vert}, \frac{1}{\vert\a+q\vert})}{G(\frac{\vert\a\rho +q\vert}{\vert\a +q\vert}, \frac{1}{\vert\a+q\vert})G(\frac{\vert q\vert}{\vert\a +q\vert} +s, \frac{1}{\vert\a+q\vert})}\cdot$$
Using (\ref{Gumbel}) and (\ref{DoubleG}) we then write
$$\log \Psi_+(s)\; =\; \kappa \,s \; +\; \int_{-\infty}^0 (e^{sx} - 1 -sx)\, (f_{1,1+\delta,\delta}+f_{a,b,\delta})(\vert x\vert)\, dx$$
for some normalizing constant $\kappa,$ with
$$ a\; =\; \frac{\vert q\vert}{\vert\a+q\vert},\quad b\; =\; \frac{\vert\a\rho+q\vert}{\vert\a+q\vert}\quad \mbox{and}\quad \delta\; =\; \frac{1}{\vert\a+q\vert}\cdot$$
We next decompose
$$(f_{1,1+\delta,\delta}+f_{a,b,\delta})(\vert x\vert)\; =\; f_{a,b+\delta,\delta}(\vert x\vert)\, +\, \frac{e^{-\vert x\vert}(1-e^{-(b-1)\vert x\vert})}{\vert x\vert (1-e^{-\vert x\vert})}$$
and deduce that $\log \Psi_+(s)$ equals
$${\tilde \kappa}\, s + \int_{-\infty}^0 (e^{sx} - 1)\frac{e^{-\vert x\vert}(1-e^{-(b-1)\vert x\vert})}{\vert x\vert (1-e^{-\vert x\vert})}\,  dx + \int_{-\infty}^0 (e^{sx} - 1 -sx)\, f_{a,b+\delta,\delta}(\vert x\vert)\, dx$$
for some normalizing constant ${\tilde \kappa}.$ Putting everything together with (\ref{IniBB}), we finally obtain the required identity
$$\AA(\a,\rho,q)\; \elaw\; \KK\;\B_{1, \frac{\a\hr}{\vert \a +q\vert}}\times\, \TT\lpa \frac{\vert q\vert}{\vert \a +q\vert}, \frac{1}{\vert \a +q\vert}, \frac{1-\a\rho}{\vert \a +q\vert}\rpa\,\times\,\TT^{-1}\lpa\frac{1-\a\hr}{\vert \a +q\vert}, \frac{1}{\vert \a +q\vert}, \frac{\a\hr}{\vert \a +q\vert}\rpa$$
with
$$\KK\; =\; \frac{\Ga(1-\a\hr)\, \times\, \Psi_+(1)}{\esp[\B_{1, \frac{\a\hr}{\vert \a +q\vert}}]}\; =\; \frac{\Ga(1-q-\a\rho)\Ga(1-\a\hr)}{\vert \a+q\vert \Ga(\vert q\vert)}\cdot$$

\qed

\begin{REMS} {\em (a) We refer to Theorems 3 and 4 in \cite{KP} for convergent series representations at infinity and asymptotic expansions at zero of the density function of $\AA(\a,\rho,q)$, obtained from the Mellin inversion formula - see \cite{K1, K2} for the case of the random variable $T.$ Beware that contrary to the fractional moments, the latter asymptotic expansion at zero does not extend to the boundary case $\{\a > 1, \rho =1-1/\a\},$ where the density is exponentially small at zero - see Formula (5.47) in \cite{PS3} for the precise asymptotic.\\
 
(b) By Formula (A.3) in \cite{LK}, it is possible to express the fractional moments of $\AA(\a,\rho,q)$ as a renormalized infinite product involving the Gamma function only. By Formula 5.3(1) p.207 in \cite{EMOT} and the Mellin inversion formula, this shows that the density function of $\AA(\a,\rho,q)$ can be interpreted as a renormalized Meijer's $G$-function having infinite parameters. This expression also entails that the Mellin transform of $\AA(\a,\rho,q)$ has poles of infinite order when $\vert \a+q\vert$ is rational, whereas all the poles are simple when $\vert\a+q\vert$ is irrational. This dichotomy was already observed and discussed in \cite{KP}.\\

(c) In the case $\delta =1,$ the formula (\ref{Malm}) dates back to Barnes - see \cite{Ba0} p.309. It is worth mentioning that this formula plays some role in non-standard limit theorems involving renormalized products of gamma random variables, in the context of random matrix theory - see Theorems 1.4 and 1.6 in \cite{NY}. Formula (\ref{Malm}) can be viewed as a L\'evy-Khintchine formula of the "third order", as observed in Remark 2.5 of \cite{NY} in the case $\delta =1,$ whereas Formula (\ref{DoubleG}) is a true L\'evy-Khintchine formula thanks to a cancellation of the quadratic term in (\ref{Malm}).
}
\end{REMS}

\section{Proof of the Corollaries}

\subsection{Proof of Corollary \ref{Doney}} Combining Parts (a) and (b)-(i) of the Theorem, (\ref{MulTT}) and the first identity in (\ref{MulTTc}), we get 
$$\frac{T}{\hat{T}}\; \elaw\; \frac{\G_\rho}{\G_{\hr}}\; \elaw \; \B_{\hr,\rho}^{-1}-\, 1.$$

\qed

\begin{REMS} {\em  Recall from Exercise VIII.3 in \cite{Be1} that in the case with negative jumps, the undershoot at time $T$ is distributed as
$$|L_T|\;\elaw\;\B_{\a\hr, 1-\a\hr}^{-1}\; -\; 1.$$
In the Cauchy case $\a=1$, we observe the interesting identities
$$\frac{T}{\hat{T}}\; \elaw\; |L_T|\; \elaw\; \SS(1,\rho)\, -\, 1\;\elaw\; \frac{1}{1-\II(1,\rho)}\, -\, 1.$$}
\end{REMS}

\medskip

\subsection{Proof of Corollary \ref{Dondon}} The argument relies on Part (b)-(iii) of the Theorem, and is the same as above. Simplifying the constants and using the well-known factorization 
\begin{equation}
\label{BeGa}
\B_{a,b}\,\times\,\G_{a+b}\; \elaw\; \G_a,
\end{equation}
we obtain
$$\frac{\AA(\a,\rho,-1)}{\AA(\a,\hr,-1)}\;\elaw \;\B_{1, \frac{\a\hr}{1- \a}}\,\times\,\G_{\frac{1-\a\rho}{1-\a}}\,\times\,\B_{1, \frac{\a\rho}{1-\a }}^{-1}\,\times\,\G_{\frac{1-\a\hr}{1-\a}}^{-1}\;\elaw \; \L \, \times\; \L^{-1}\;\elaw\; \U^{-1} - 1.$$
\qed

\subsection{Proof of Corollary \ref{Subutex}} This follows from Part (b)-(iii) of the Theorem at the limit $\{\a < 1, \hr =0\},$ recalling the conventions $\B_{1,0}\equiv\TT_{a,b,0}\equiv 1.$ 

\qed

\begin{REMS} {\em (a) Letting $\a\to 1,$ we recover
$$\int_0^\infty \frac{dt}{(1+t)^q}\; =\; \frac{1}{q-1}\cdot$$
Letting $\a\to 0,$ we obtain
\begin{equation}
\label{Cresson}
\int_0^\infty \frac{dt}{{(1+\sigma^{(\a)}_t)}^q}\; \claw\;\L
\end{equation}
for every $q >0,$ in accordance with the fact that $\sigma^{(0)}$ can be interpreted as a pure murder at an exponential time $\L.$ 

\medskip

(b) Letting $q\to+\infty$ and (\ref{MelTT}) and (\ref{Gumbel}), one can check that 
\begin{equation}
\label{Gumgum}
\int_0^\infty e^{-\sigma^{(\a)}_t}\,dt\; \elaw\;\L^{[1-\a]}
\end{equation}
where $\{\log \L^{[t]},\, t\ge 0\}$ is the real L\'evy process starting from zero characterized by 
$$\log \L^{[1]}\; \elaw \;\log\L,$$ 
and which for this reason might be called the Gumbel-L\'evy process, recalling that $\log\L$ follows a standard Gumbel distribution. Notice that (\ref{Gumgum}) can also be obtained from Part (b)-(i) of the Theorem with $\{\a < 1, \rho =0\}$ and $q\to +\infty.$ Using the independence and stationarity of the increments of the Gumbel-L\'evy process, we retrieve the following factorizations of the standard exponential law
$$\lpa\int_0^\infty e^{-\sigma^{(\a)}_t}\,dt\rpa\times\lpa\int_0^\infty e^{-\sigma^{(1-\a)}_t}\,dt\rpa\; \elaw\;\L,$$
which had been observed in Formula (14) of \cite{BY1}. It is interesting to compare these factorizations with those for the standard Pareto law obtained in Corollary 2. It is also easy to deduce the identity (\ref{Gumgum}) from the results of \cite{BY1}.\\

(c) Taking $q=1$, we notice that the Riemann random integral has a Weibull distribution:
\begin{equation}
\label{Weibull}
\int_0^\infty \frac{dt}{1+\sigma^{(\a)}_t}\; \elaw\; \frac{1}{1-\a}\; \L^{1-\a}.
\end{equation}
This can also be deduced from the results of \cite{Be2} - see Lemma 2 (ii), Corollary 3 (ii) and Comment 4 therein. We then retrieve (\ref{Cresson}) for $q=1$ in letting $\a\to 0.$ This also shows that
$$(1-\a)\int_0^\infty \frac{dt}{1+\sigma^{(\a)}_t}\; \claw\; 1, \qquad \a\to 1.$$
We last deduce from (\ref{Weibull}), Shanbhag-Sreehari's identity - see e.g. Exercise 29.16 in \cite{S}, and Corollary 2 in the case $\rho =0,$ the identity
\begin{equation}
\label{ShSr}
\int_0^T \frac{dt}{1-\sigma^{(\a)}_t}\; \elaw\; \frac{1}{1-\a}\; \L\, \times\,{\Z}_{1-\a}^{1-\a}.
\end{equation}
This can also be derived directly from Part (b)-(iii) of the Theorem.\\

(d) In \cite{TW}, it is observed that the law of the above Riemann random integral is a size-bias of an inverse generalized stable distribution, that is its density is the unique solution to an integro-differential equation of first order involving a Beta random variable.}

\end{REMS}

\medskip

\subsection{Proof of Corollary \ref{dual}} To show Part (a), we apply first the factorization (\ref{MulTT}) to obtain
$$\TT\lpa 1, \frac{1}{\a +q}, \frac{1-\a\hr}{\a +q}\rpa\;\elaw\; \TT\lpa 1, \frac{1}{\a +q}, \frac{1-\a\hr'}{\a +q}\rpa \,\times\,\TT\lpa\frac{1+q+\a\rho'}{\a +q}, \frac{1}{\a +q}, \frac{\a(\rho-\rho')}{\a +q}\rpa$$
and
$$\TT\lpa \frac{\a\hr}{\a +q}, \frac{1}{\a +q}, \frac{\a\rho}{\a +q}\rpa\;\elaw\; \TT\lpa \frac{\a\hr'}{\a +q}, \frac{1}{\a +q}, \frac{\a\rho'}{\a +q}\rpa \,\times\,\TT\lpa\frac{\a\hr}{\a +q}, \frac{1}{\a +q}, \frac{\a(\rho-\rho')}{\a +q}\rpa.$$
The conclusion follows then from the Theorem. To obtain Part (b), we use an analogous argument together with the following consequence of (\ref{MulBB}):
$$\B_{1, \frac{\a\hr}{\vert \a +q\vert}}\;\elaw\; \B_{1, \frac{\a\hr'}{\vert \a +q\vert}} \,\times\,\B_{\frac{\vert\a\rho' +q\vert}{\vert \a +q\vert},\frac{\a(\rho'-\rho)}{\vert \a +q\vert}}.$$ 
We omit the details.

\qed

\begin{REMS} {\em (a) Taking $\{\a > 1, \rho' = 1-1/\a, q = 0\}$ entails that for every $\a > 1,$ one has the factorization
\begin{equation}
\label{GJ}
T \; \elaw\; \XX\,\times\, \Z_{\frac{1}{\a}}.
\end{equation}
with
$$\XX\; =\; \lpa\frac{\Ga(1+\a\rho)\Ga(\a\hr)\Ga(1+1/\a)}{\Ga(1+\a)\Ga(\a)}\rpa\times\;\TT\lpa 1, \frac{1}{\a}, \frac{1}{\a} -\hr\rpa\;\times\; \TT^{-1}\lpa \hr, \frac{1}{\a}, \frac{1}{\a}-\hr\rpa.$$
This had been observed in Corollary 4 of \cite{GJ}, with a different expression for the factor $\XX$ involving Darling's integral. More precisely, the density function of $\XX$ in \cite{GJ} reads 
\begin{equation}
\label{Darling}
\frac{\sin(\pi\a\hr)\,g_{\a,\rho}(x^{\frac{1}{\a}})}{\pi\a\, (x^2 +2x \cos(\pi\a\hr) +1)}
\end{equation}
with the notation
$$\log g_{\a,\rho}(x)\; =\; \frac{\sin(\pi\rho)}{\pi}\,\int_0^\infty\frac{x\log(1+ u^{\frac{1}{\a}})}{u^2 +2xu \cos(\pi\rho) +1}\, du.$$ 
However, it seems difficult to deduce from this expression of its density the factorization of $\XX$ as the above double infinite Beta product.\\

(b) In the spectrally negative case $\rho=1/\a$, the density of $\XX$ in (\ref{GJ}) has the simple expression
$$ \frac{(-\sin \pi\a)(1+x^{1/\a})}{\pi\a(x^2 - 2x \cos \pi\a  +1)}\cdot$$ 
This had been observed in \cite{TS0} - see Theorem 3 therein, thanks to a computation on the classical Mittag-Leffler function. We will come back to this example in Section 4.\\

(c) Using the same argument as in the proof of Corollary \ref{dual}, one can prove the following extension of (\ref{GJ}):  there exists an explicit infinite double Beta product $\XX$ such that
$$\AA(\a,\rho,q)\; \elaw\; \XX\,\times\, \Z_{\frac{1}{\a +q}}$$
as soon as $1 \le q + \a \le 1 +\a\rho.$} 

\end{REMS}

\subsection{Proof of Corollary \ref{WH}} In the spectrally positive case, the factorization is obvious by (\ref{ExpSP}), taking $Z =(\a+q)\xi$ and $\sigma_t =t.$ Assume next that $L$ is not spectrally positive. If $q =-\a,$ the factorization is immediate from Part (b)-(ii) of the Theorem, taking $Z\equiv 0$ and for $\sigma$ a pure murder at the exponential time $r\L.$ \\

Suppose now $q > -\a.$ Using the notation of Theorem 6.2 (b)-(i) and applying the concatenation formula given as Formula (25) in \cite{KP}, we see that for every $s>-\a/(\a+q)$:
\begin{eqnarray*}
\frac{s\Psi_-(-(s+1))}{\Psi_-(-s)} & = & \frac{s \,G(\frac{\a}{\a +q} +s+1,\frac{1}{\a +q})G(\frac{\a\hr}{\a +q} +s, \frac{1}{\a +q})}{G(\frac{\a}{\a +q} +s,\frac{1}{\a +q})G(\frac{\a\hr}{\a +q} +s+1,\frac{1}{\a +q})}\\
& = & \frac{s\,\Ga(\a +(\a +q)s)}{\Ga(\a\hr + (\a +q)s)}\\
& = & \frac{\Ga(\a)\, s}{\Ga(\a\hr)} \, +\, \int_{-\infty}^0 \!(e^{s x} - 1 - sx)\frac{\a\rho\, e^{\frac{\a x}{\a +q} }(\a+(1-\a\hr)e^{\frac{x}{\a +q}})}{(\a+q)^2\Ga(1-\a\rho)(1-e^{\frac{x}{\a +q}})^{\a\rho+2}}\, dx
\end{eqnarray*}
where in the third equality we have used the Lemma in \cite{BS1} - see also example 2 p. 1688 in \cite{KR}. This implies
$$\frac{s\Psi_-(-(s+1))}{\Psi_-(-s)}\; =\; \log\esp[e^{s Z_1}],$$
where $\{Z_t, \, t\ge 0\}$ is a spectrally negative L\'evy process with explicit characteristics, positive mean and without Gaussian part. Applying Bertoin and Yor's criterion - see Proposition 2 in \cite{BY2} - entails
$$\Psi_-(s)\; =\; \esp\lcr \lpa \int_0^\infty e^{-Z_t}\, dt\rpa^s\,\rcr.$$
In order to identify the factor $\Psi_+(s),$ we use again Formula (25) in \cite{KP} to  obtain
\begin{eqnarray*}
\frac{s\Psi_+(s-1)}{\Psi_+(s)} & = & \frac{G(\frac{1-\a\hr}{\a +q} +s,\frac{1}{\a +q})G(\frac{1}{\a +q} +s+1, \frac{1}{\a +q})}{G(\frac{1-\a\hr}{\a +q} +s+1,\frac{1}{\a +q})G(\frac{1}{\a +q} +s,\frac{1}{\a +q})}\\
& = & \frac{\Ga(1 +(\a+q)s)}{\Ga((1-\a\hr) +(\a+q)s)}\\
& = & \frac{1}{\Ga(1-\a\hr)}\lpa 1 \, +\, \int_0^\infty \!(1-e^{-s x})\,\frac{\a\hr\, e^{-\frac{x}{\a +q}}}{(\a+q)(1-e^{\frac{x}{\a +q}})^{1+\a\hr}}\, dx\rpa.
\end{eqnarray*}
This implies that for every $s>0,$ one has
$$\frac{s\Psi_+(s-1)}{\Psi_+(s)}\; =\; -\log\esp[e^{-s \sigma_1}],$$
where $\{\sigma_t, \, t\ge 0\}$ is a subordinator with explicit characteristics, murdered but without drift. Applying Carmona-Petit-Yor's criterion - see Proposition 3.3 dans \cite{CPY} - shows that
$$\Psi_+(s)\; =\; \esp\lcr \lpa \int_0^\infty e^{-\sigma_t}\, dt\rpa^s\,\rcr.$$
Putting everything together entails the required factorization
$$\AA(\a,\rho,q)\; \elaw\; \lpa \int_0^\infty e^{-Z_t}\, dt\rpa \times \lpa \int_0^\infty e^{-\sigma_t}\, dt\rpa.$$
The case $q < -\a$ is entirely analogous and we omit the details. The L\'evy-Khintchine exponents of $Z$ and $\sigma$ are here respectively given by
$$\frac{s}{\Ga(1-\a\hr)} \, +\, \int_{-\infty}^0 \!(e^{s x} - 1 - sx)\frac{\a\hr\, e^{\frac{x}{\vert \a +q\vert} }(1+ \a\hr\, e^{\frac{x}{\vert \a +q\vert}})}{\vert \a+q\vert^2\Ga(1-\a\hr){(1-e^{\frac{x}{\vert \a +q\vert}})}^{\a\hr+2}}\, dx$$
and
\begin{equation}
\label{Riem}
\frac{\Ga(\a)}{\Ga(\a\hr)}\; +\; \frac{\a\rho}{\vert \a +q\vert \Ga(1-\a\rho)} \,\int_0^\infty \!(1-e^{-s x})\,\frac{e^{-\frac{\a x}{\vert \a +q\vert}}}{{(1-e^{\frac{x}{\vert \a +q\vert}})}^{1+\a\rho}}\, dx.
\end{equation}
\qed

\begin{REMS} {\em From the above proof, it is interesting to observe that the Riemann random integral of Corollary \ref{Subutex} can also be viewed as the perpetuity of a subordinator. One has
$$\int_0^\infty \frac{dt}{{(1+\sigma^{(\a)}_t)}_{}^q}\;\elaw\;  \int_0^\infty e^{-\sigma^{(\a,q)}_t}\, dt$$
where $\{\sigma^{(\a,q)}_t, \, t\ge 0\}$ is the (non-killed) subordinator whose L\'evy-Khintchine exponents are given by (\ref{Riem}) with $\hr =0.$ Observe that
$$\{\sigma^{(\a,q)}_{q^{-\a}t}, \, t\ge 0\}\;\claw\; \{\sigma^{(\a)}_t, \, t\ge 0\}\qquad\mbox{as $q\to +\infty,$}$$
as expected.}

\end{REMS}

\smallskip

\subsection{Proof of Corollary \ref{AD}} It is well-known and easy to see from (\ref{IniBB}) that $\log \B_{a,b}$ is ID for all $a,b >0.$ Since infinite divisibility is preserved under weak convergence, we see from the Theorem that $\log \AA(\a,\rho,q)$ is also ID. To obtain its self-decomposability, it is sufficient to prove the property for the factors $\log\TT(a,b,c)$ and $\log\B_{a,b}$ involved in the  factorization. It is known and not difficult to see that  
$$\log \B_{a,b}\;\, \mbox{is SD}\quad\Leftrightarrow\quad x\mapsto \frac{x^a(1-x^b)}{1-x}\;\;\mbox{is non-decreasing on $(0,1)$}\quad\Leftrightarrow\quad 2a+b\geq 1.$$
On the other hand, it is clear from (\ref{MelTT}) that for every $a,b,c >0,$ one has 
\begin{equation}
\label{A&M}
x\mapsto \frac{x^a(1-x^c)}{(1-x)(1-x^b)}\;\;\mbox{is non-decreasing on $(0,1)$}\quad\Leftrightarrow\quad \log \TT(a,b,c)\;\, \mbox{is SD,}
\end{equation}
and that the condition on the left-hand side is fulfilled as soon as $2a +c \ge 1\wedge b.$ Combining these facts, it is easy to deduce from the Theorem that $\log \AA(\a,\rho,q)$ is SD as soon as $\a(2-\rho)\ge (\a+q)\wedge 1.$ On the other hand, one can show from (\ref{A&M}) after some computations that $\log \TT(1/4,1,1/4)$ is not SD, so that  $\log\AA(1/2,1/2,1/2)$ is not SD either, by Part (b)-(i) of the Theorem.

\qed

\begin{REMS}{\em (a) It does not seem that the criterion on the left-hand side of (\ref{A&M}) can be conveniently characterized in terms of $(a,b,c),$ which would allow to characterize the self-decomposability of $\log\AA(\a,\rho,q)$ in terms of $(\a,\rho,q).$ See Theorem 1 and Remark 2 in \cite{AM} for a similar question related to the self-decomposability of higher order of the random variable 
$\log \G_t.$\\

(b)  The multiplicative infinite divisibility and self-decomposability of the perpetuity associated to a subordinator have been characterized in \cite{HY} - see Proposition 3.3 and Proposition 3.5 therein. It does not seem that there exists any criterion in the literature for the multiplicative self-decomposability of the exponential functional of a general L\'evy process.}
\end{REMS}

\medskip

\subsection{Proof of Corollary \ref{Chinchin}} First, it is clear from Part (b)-(ii) of the 
Theorem that $\AA(\a,\rho,q)$ has a decreasing density which is bounded at zero in the case $q=-\a.$ \\

Suppose now $1-\a\hr \ge\a + q > 0.$ We see from Part (b)-(i) of the Theorem that $\AA(\a,\rho,q)$ admits the random variable 
$$\B_{1, \frac{1-\a\hr}{\a +q}}\; \elaw\; \U\,\times\,\B_{2,\frac{1-\a\hr}{\a +q}-1}$$
as a multiplicative factor, where the identity in law stems from (\ref{MulBB}). By Khintchine's theorem - see e.g. Exercise 29.21 in \cite{S}, this entails that $\AA(\a,\rho,q)$ has a non-increasing density. Moreover, it is clear that in this situation the first negative pole of $s\mapsto \esp[\AA(\a,\rho,q)^s]$ is then at $s =-1,$ and is  simple. An application of the residue theorem entails then that the density function of $\AA(\a,\rho,q)$ is necessarily bounded at zero. \\

Suppose next $\a + q >1-\a\hr \ge 0.$ From the proof of Part (b)-(i), we deduce the factorization 
$$(1+s)\esp[\AA(\a,\rho,q)^s] \; = \; \frac{\Ga(2+s)}{\Ga(1+ \frac{1-\a\hr}{\a+q} +s)}\,\times\,\Psi(s),$$
where $\Psi$ is an analytic function on $(-(2\wedge(1+1/(\a+q))); \a\hr/(\a +q)).$ This shows that the function $s\mapsto (1+s)\esp[\AA(\a,\rho,q)^s]$ is analytic on $(-(2\wedge(1+1/(\a+q))); \a\hr/(\a +q))$ and vanishes on
$$s\; = \;-1-\frac{ (1-\a\hr)}{(\a+q)}\;\in\;(-(2\wedge(1+1/(\a+q))), \a\hr/(\a +q)).$$ 
Hence, this function cannot be the Mellin transform of any positive random variable, which means that $\AA(\a,\rho,q)$ does not admit $\U$ has a multiplicative factor. Using the only if part of Khintchine's theorem, this implies that its density is not non-increasing. The criterion for the case $\{\rho\neq 0, \a +q<0\}$ is proved in an entirely analogous way, and we omit the details. In the case $\{\rho = 0, \a +q<0\},$ we obtain  from (\ref{BeGa}) the factorization
\begin{equation}
\label{Mixt}
\AA(\a,0,q)\;\elaw\; \Ga(1-\a)\,\L\;\times\; \TT^{-1}\lpa \frac{1-\a}{\vert \a +q\vert}, \frac{1}{\vert \a +q\vert}, \frac{\a}{\vert \a +q\vert}\rpa,
\end{equation}
which shows that $\AA(\a,0,q)$ is a mixture of exponentials and hence has a completely monotone density (see Theorem 51.6 and Proposition 51.8 in \cite{S}).
\qed

\subsection{Proof of Corollary \ref{AGGC}.} By Part (a) of the Theorem, we have the identity
$$\AA(\a,1-1/\a,q)\;\elaw\;\frac{1}{(\a +q) \Ga(\a)}\times\;\prod_{n=0}^\infty \lpa \frac{n+\a}{n+1}\rpa \B^{-1}_{\frac{n+1}{\a+q},\frac{\a-1}{\a+q}}.$$
By Theorem 2.2 in \cite{BS2}, we know that all factors on the right-hand side have GGC distributions, and we can conclude by the main result of \cite{BJTP} and the stability of the GGC property at the weak and non-degenerate limit.

\qed

\begin{REMS} {\em (a) The above argument shows that $\AA(\a,\rho,q)$ is always the independent quotient of two random variables having a GGC distribution. However, such an independent quotient is not necessarily a GGC. For example, one can show from Theorem 4.1.1 in \cite{Bd} that the quotient of two uniform random variables, whose density is constant and equal to $1/2$ in $(0,1),$ cannot be a GGC.\\

(b) The above argument also shows that the random variables $\AA(\a,0,q)^{-1}, q > -\a$ and $\AA(\a,1,q)^{-1}, q < -\a$ have a GGC distribution. Taking the limit $q\to +\infty$ resp. $q\to -\infty$, we deduce that the law of the reciprocal of the perpetuity
$$\int_0^\infty e^{-\sigma^{(\a)}_t}\,dt,$$
with the notation of Corollary \ref{Subutex}, is also a GGC. This property can also be obtained via the method of Theorem 4 in \cite{BS2}. It is worth mentioning that on the other hand, the above perpetuity is never infinitely divisible, except in the degenerate case $\a =1$. Indeed, it follows from (\ref{Gumgum}) and Stirling's formula that the moment generating function of its $t$-th power reads
$$\sum_{n\ge 0} \frac{\Ga(1+tn)^{1-\a}}{n!} \,\lbd^n$$
and is finite for all $\lbd\in\rl$ and $t < 1/(1-\a),$ so that the upper tails of its distribution function are superexponentially small.  \\

(c) The GGC random variable $\AA(\a,1-\a,q)$ has finite negative moments of any order. By Theorem 4.1.4 in \cite{Bd}, this entails that its Thorin measure has infinite mass, with the terminology of \cite{Bd}. Using the full extent of Theorem 2.2 in \cite{BS2}, one can also show that all positive powers $\AA(\a,1-1/\a,q)^s$ with magnitude $s \ge \sup (1/2, q/(\a+q))$ have a GGC distribution with infinite Thorin measure.}
\end{REMS}

\section{Further remarks}

In this last section we revisit known results and we display some new particular cases where the law of $\AA(\a,\rho,q)$ has a simpler expression. We also consider the hitting times for points in the relevant situation $\a >1,$ complementing the results we had obtained in our previous study \cite{LS1}.

\subsection{The stopped extrema} It is easy to deduce from (\ref{IniBB}) et (\ref{MelTT}) the convergences in law
$$\TT(a\eps,\eps,c\eps)^\eps\; \claw\; \B_{a,c}\qquad\mbox{and}\qquad \TT(a,\eps,c\eps)^\eps\;\claw\;1$$
when $\eps\to 0,$ for every $a,c >0.$ Using the Theorem, we can then recover the identities (\ref{Rog}) and (\ref{Port}) on stopped extrema, which we recall to be read as
$$\sup\{L_t, \, t<T\}\;\elaw\;\B_{\a\hr, \a \rho}^{-1},\qquad\mbox{and}\qquad\inf\{L_t, \, t<T\}\;\elaw\;\B_{1-\a\hr, \a\hr}.$$
It is worth mentioning that $\sup\{L_t, \, t\le T\} = \sup\{L_t, \, t <T\},$ whereas
$$\inf\{L_t, \, t\le T\} \; =\; L_T\; \elaw\; 1\, -\, \B_{\a\hr, 1-\a \hr}^{-1}$$
(see the above Remark 6). It would be quite interesting to investigate the law of the stopped functionals
$$\sup\{\vert L_t -1\vert, \, t<T\}\qquad\mbox{and}\qquad\sup\{L_t -L_s, \, s, t<T\},$$
which describe respectively the two-sided supremum and the range of the process until the stopping time $T.$ It is well-known among specialists that the law of the corresponding unstopped functionals is very complicated to study in general even at the asymptotic level - see however \cite{F} for the Brownian case.

\subsection{The Cauchy case} We have here  $\{\a =1, \rho\in (0,1)\}$ and we set $L = \{C^{(\rho)}_t,\, t\ge 0\}$ for the corresponding Cauchy process starting from one. In the symmetric case, the density of the first passage time $\AA(1, 1/2,0)$ has been expressed by Darling as a certain running integral - see the last paragraph of \cite{Da}. It should be possible to recover this formula from our factorization, and also to extend it in the asymmetric case, but we have not investigated this question.\\

In the case $q=-\rho,$ we deduce from Part (b) (i) of the Theorem and the second  identity in (\ref{MulTTd}) the identity.
$$\AA(1, \rho,-\rho)\;\elaw\; \frac{1}{\hr} \lpa \frac{\Z_{\hr}}{\Z_{\hr}}\rpa^{\hr}.$$
Thanks to a standard computation - see Exercice 4.21 (3) in \cite{CY}, this identity can be restated with the help a cut-off dual Cauchy random variable, showing that  $\AA(1, \rho,-\rho)$ has an explicit density: one has
$$\int_0^T \frac{\hr\,dt}{{\vert C^{(\rho)}_t\vert}_{}^\rho}\; \elaw\; (C^{(\hr)}_1-1)\,\vert\, C^{(\hr)}_1 >1\; \sim\;\frac{\sin(\pi \hr)}{\pi \hr (x^2+2x \cos(\pi \hr)+1)}\cdot$$
In particular, this entails from the results in \cite{BHC} that $\AA(1, \rho,-\rho)$ is ID as soon as $\rho \ge 1/2,$ and from those of \cite{Di} that $\AA(1, 1/2,-1/2)$ is SD. Observe also that
$$\hr\,\AA(1, \rho,-\rho)\;\claw\; \frac{\L}{\L}\; \sim\;\frac{1}{(x+1)^2}$$
as $\rho \to 1,$ whereas $\AA(1, \rho,-\rho)\claw 1$ when $\rho \to 0,$ as expected. In the case $q=-\hr,$ we have the analogous identity
$$\AA(1, \rho,-\hr)\;\elaw\; \frac{1}{\rho}\lpa \frac{\G_{\rho}}{\G_{\hr}}\rpa^{\rho},$$
to be compared with Corollary 1. The infinite divisibility of the random variable on the right-hand side is an open question when $\rho\neq 1/2.$ \\

We finally observe from the proof of the Theorem the identity $\xi^{(\rho)} \elaw -\xi^{(\hr)}$ between the hypergeometric L\'evy processes respectively associated to $C^{(\rho)}$ and $C^{(\hr)}$. This implies
$$\AA(1, \rho,q)\;\elaw\;\AA(1, \hr, -2-q).$$
In particular, we have the identities
$$\int_0^T \frac{\rho\,dt}{{\vert C^{(\rho)}_t\vert}_{}^{\rho+1}}\; \elaw\; (C^{(\rho)}_1-1)\,\vert\, C^{(\rho)}_1 >1\qquad\mbox{and}\qquad \int_0^T{\vert C^{(\rho)}_t\vert}_{}^{-2} dt\; \elaw\; {\hat T}$$
for every $\rho \in (0,1).$ In the limiting cases we find
$$\rho\,\AA(1, \rho,-\rho-1)\;\claw\; \frac{\L}{\L}\quad\mbox{as $\rho \to 0$} \qquad\mbox{and}\qquad \int_0^\infty \frac{dt}{(1+t)^2}\; =\; 1.$$

\subsection{Some further explicit identities} In this paragraph, we describe some other situations where the law of $\AA(\a,\rho,q)$ has a relatively simple expression - see (\ref{ML1}), (\ref{Z2}) and (\ref{C3}) below. We restrict ourselves to the cases where $q$ is a non-negative integer and $-L$ is a subordinator or $L$ is spectrally one-sided, although other simplifications are certainly possible. We will state our results without proof, for the sake of conciseness. These proofs are simple rearrangements relying on (\ref{Frac}), the Legendre-Gauss multiplication formula for the Gamma function, concatenation formul\ae\, for the double Gamma function to be found in (A1) and (A8) of \cite{LK}, and fractional moment identifications. They have been typesetted and are available upon request. 

\subsubsection{The case where $-L$ is a subordinator} We suppose here $\{\a\in (0,1), \rho =0\}.$ Introduce the so-called Mittag-Leffler random variable
$$\M_\a\; =\; \Z_\a^{-\a},$$
whose denomination comes from the fact that the Laplace transform of $\M_\a$ is expressed in terms of the classical Mittag-Leffler function - see Exercise 29.18 in \cite{S}. It is well-known, easy to see and mentioned in the introduction that 
$$\AA(\a,0,0)\;\elaw\;\M_\a.$$ 
We can prove the following generalization, for every $q\in\NN:$ 
\begin{equation}
\label{ML1}
\AA(\a, 0, q)\;\elaw\; (\a +q)^q\,\times\, \M_{\frac{\a+q}{1+q}}\,\times\,\cdots\,\times\,\M_{\frac{\a+q}{1+q}}^{(\frac{q}{\a+q})}.
\end{equation}
Recall that in the case $q\le -\a$ we have also the factorization (\ref{Mixt}). The latter does not seem to take a simpler form except in the case $q=-1$ - see (\ref{ShSr}) above. However, its interesting feature is the GGC property for $\AA(\a, 0, q),$ which is then the product of two GGC random variables. 

\subsubsection{The case where $L$ is spectrally positive} We assume here $\{\a >1, \rho =1-1/\a\}.$ It is well-known, easy to see and mentioned in the introduction that 
$$\AA(\a,1-1/\a,0)\;\elaw\;\Z_{\frac{1}{\a}}.$$ 
We can prove the following generalization, for every $q\in\NN:$ 
\begin{equation}
\label{Z2}
\AA(\a, 1-1/\a, q)\;\elaw\; (\a +q)^q\,\times\, \Z_{\frac{1+q}{\a+q}}\,\times\,\cdots\,\times\,\Z_{\frac{1+q}{\a+q}}^{(\frac{q}{\a+q})}.
\end{equation}

\subsubsection{The case where $L$ is spectrally negative} We suppose here $\{\a >1, \rho =1/\a\}$ 
For every $\mu\in (0,1),$ introduce the $\mu-$Cauchy random variable $\C_\mu$ with density function
$$\frac{\sin(\pi\mu)}{\pi\mu (x^2 + 2 \cos (\pi\mu) x + 1)}.$$
Recall that with the notation of Section 4.2, one has $\C_\mu\elaw C^{(\mu)}\,\vert\, C^{(\mu)} > 0.$ In the context of Corollary \ref{dual} we can prove the following refinement of the main result of \cite{TS0} - see also Remark 7 (b) above. For every $q\in\NN,$ one has
\begin{equation}
\label{C3}
\AA(\a, 1/\a, q)\;\elaw\; \lpa \C_{\frac{\a+q}{2+q}}\,\times\,\cdots\,\times\,\C_{\frac{\a+q}{2+q}}^{(\frac{1+q}{\a+q})}\rpa\times\,\AA(\a, 1-1/\a, q).
\end{equation}
As already mentioned, in the case $q=0$ the density of $\C_{\frac{\a}{2}}\times\C_{\frac{\a}{2}}^{(\frac{1}{\a})}$ is explicit and reads
$$\frac{(-\sin \pi\a)(1+x^{1/\a})}{\pi\a(x^2 - 2x \cos \pi\a  +1)}\cdot$$
However, we could not derive in general a simple expression of this kind for the density of the independent factor
$$\C_{\frac{\a+q}{2+q}}\,\times\,\cdots\,\times\,\C_{\frac{\a+q}{2+q}}^{(\frac{1+q}{\a+q})}.$$

\subsection{Hitting times} In this last paragraph we consider the first hitting time of zero for $L:$
$$\tau(\a,\rho)\; =\; \inf\{t >0, \; L_t =0\},$$
which is finite a.s. if and only if $\a > 1$ - see e.g. Example 43.22 in \cite{S}. In \cite{LS1}, Formula (8), we have derived the following factorization
\begin{equation}
\label{Hit}
\tau(\a,\rho)\; \elaw\; \C_{\a\hr}^{(\frac{1}{\a})}\,\times\; \Z_{\frac{1}{\a}}.
\end{equation}
Observe that the above $\C_{\a\hr}^{(\frac{1}{\a})}$ has explicit density
$$\frac{\sin (\pi\a\hr)\sin(\frac{\pi}{\a})\,x^{\frac{1}{\a}}}{\pi \sin(\pi\hr)(x^2 +2 x\cos(\pi\a\hr) + 1)},$$
which should be compared with the more complicated density (\ref{Darling}) involved in the factorization of $T(\a,\rho)$. In \cite{LS1}, the factorization (\ref{Hit}) was obtained from the main result of \cite{KKMW} computing the fractional moments of $\tau$ via the perpetuity approach. It was also shown that this computation can be quickly derived from the potential formula - see Section 2.4 in \cite{LS1}. \\

Using Proposition \ref{GZT} (b), the simple factorization (\ref{Hit}) can be rewritten as
$$\frac{1}{\Ga(1+\a)}\;\TT\lpa 1+\frac{1}{\a}, \frac{1}{\a\hr}, \frac{1}{\a\hr} -1\rpa\;\times\;\TT^{-1}\lpa 1+\frac{1}{\a}, \frac{1}{\a\hr}, \frac{1}{\a\hr} -1\rpa\;\times\;\TT^{-1}\lpa \frac{1}{\a}, \frac{1}{\a}, 1-\frac{1}{\a}\rpa.$$
This "infinite triple Beta product" looks more complicated than the factorizations we have obtained for $\AA(\a,\rho,q).$ It also allows to prove the following analogue of Corollary \ref{AD}, which we state together with a factorization generalizing the main result of \cite{TS1}, where it was established in the case $\rho = 1/\a = 1-\rho'.$ We omit the proofs.

\begin{PROP} {\em (a)} For every $\a > 1$ and $\rho > \rho',$ one has
$$\tau (\a,\rho)\; \elaw\; {(\C_{\frac{\hr}{\hr'}}^{\a\hr'})}^{(\frac{1}{\a})}\times\; \tau (\a, \rho').$$

{\em (b)} The random variable $\log \tau(\a,\rho)$ is {\em ID.} Moreover, it is {\em SD} whenever $\rho \ge 1/2.$
\end{PROP}

It is probably true that $\log \tau(\a,\rho)$ is always SD, but the argument to prove this involves the sum of two functions of the type (\ref{A&M}), and eludes us. We also believe that $\tau(\a,\rho)$ {\em itself} is ID. In the case $\rho =1/2,$ using the first formula (5.12) in \cite{Y3} and Bochner's subordination, this would entail that the real $\a-$Cauchy distribution ($1<\a\le 2$), with explicit density function
$$\frac{\a\sin(\pi/\a)}{2\pi (1+\vert x\vert^\a)}$$
over $\rl,$ is also ID. This interesting question is stated as an open problem in Remark 2.9 (i) of \cite{Y3}.
 
\bigskip

\noindent
{\bf Acknowledgements.} We thank Jean Bertoin and Lo\"\i c Chaumont for references.

\end{document}